\documentclass[12pt,a4paper,reqno]{amsart}

\usepackage{graphicx}   
\usepackage{graphics}
\usepackage{epsfig}
\usepackage{anysize}
\marginsize{2.5cm}{2.5cm}{2.5cm}{2.5cm}

\begin{document}
\title[Niftiyev A.A.: Correctness of the optimal control problems for ...]{Correctness of the optimal control problems for
distributed parameter systems \\
(\footnotesize\textit{survey})}
\author{Niftiyev A.A.}
\maketitle

\bigskip

\begin{center}
\textit{Institute of Applied Mathematics Baku State University}
\end{center}

\begin{center}
\textit{Khalilov str.23, Az1148 Baku Azerbaijan}
\end{center}

\begin{center}
\textit{\underline {aniftiyev@yahoo.com}}
\end{center}

\bigskip

The problem of existence of optimal control for nonlinear processes, in
contradistinction to linear, is investigated a little. This problem, for the
processes, described by the ordinary differential equations, is studied, as
a rule under Philipov condition ([1]) - at condition of convexity of
admissible speeds set. However, it is known that, this condition covers very
narrow class of nonlinear problems, even in some cases of "almost linear"
problems. For multivariate variational problems and the optimal control
problems in the processes described by the equation with the partial
derivatives, in this direction we note works [2-7].

The basic idea of these works consists of providing weak semi-continuity of
the functional on weekly compact set.

The second direction in the existence of the solution is the proof of
`individual theorems'' of the existence which take into account specificity
of particularly considered problem. In this direction are received the most
important results in the works [8-10].

The essence of this direction is connected to necessary conditions of an
optimality.

The problem of existence of optimal control is investigated also in [11-13].

In the classes of problems considered in this works, the functional depends
on parameter and existence of optimal control is proved for values of
parameter from dense set . Though such type results have been obtained at
enough general conditions, the fulfillment of these conditions for
particularly taken parameter, generally speaking, is difficult.

Here an original idea is suggested to prove the existence of optimal control
for some types of non- linear problems. The obtained results can be
considered as individual existence theorems (in some sense).

The idea of the proof consists of the following:

Any minimizing sequence is chosen. It's clear, that this sequence is
minimizing for indefinitely many other functionals. If there is strict
convex among these functionals it turns out strong convergence of this
sequence. (exactly the subsequence). After it under natural conditions
theorems of existence are proved. The properties of minimizing sequence are
used for finding such strict convex functionals.

The following problems have been considered

\bigskip

1. Non convex variational problems;

2. The optimal control problems described by the first and second order
non-linear equations with abstract evolution;

3. The optimal control problem for the system described by the Gursa-Darbu
equation;

4. The optimal control problem for the system described by the parabolic
type equation.

\bigskip

\section*{At first considered non convex variational problems [16]}

\begin{equation}
\label{eq1}
J\left( {x} \right) = \int\limits_{a}^{b} {f\left( {t,x\left( {t}
\right),\dot {x}\left( {t} \right)} \right)dt} \to min,
\end{equation}

\begin{equation}
\label{eq2}
x_{i} \left( {c_{i}}  \right) = d_{i} ,i \in I.
\end{equation}

Here $a,b$ are the given numbers

\[
t \in \left[ {a,b} \right]
x(t) = (x_1 (t),x_2 (t),...,x_n (t)),
\]

\[
\dot {x}\left( {t} \right) = \left( {\dot {x}_{1} \left( {t} \right),\dot
{x}_{2} \left( {t} \right),...,\dot {x}_{n} \left( {t} \right)} \right)
d_{i} \in R
i \in I \subseteq \left\{ {1,2,...,n} \right\},
\]

 $c_{i} $ accept values $a$
or $b$, $f\left( {t,x,s} \right)$ is the given function from $\left[ {a,b}
\right] \times R^{n} \times R^{n}$ to R, $p \ge 2.$

The set of indexes\textit{}  $i$\textit{} is defined as $I_{0} $
at which

\[
f\left( {t,x_{1} ,x_{2} ,...,x_{n} ,\dot {x}_{1} ,\dot {x}_{2} ,...,\dot
{x}_{n}}  \right)
\]

\noindent doesn't depend on a component\textit{} $x_{i}
$\textit{.} Let $I_{0} \neq \emptyset$ the following conditions
also are fulfilled:

1) $J_{\ast}  = infJ\left( {x} \right) > - \infty $;

\noindent
and $J\left( {x} \right) \to + \infty $\textit{,}  if $\left\| {x}
\right\|_{W_{p}^{1} \left( {a,b} \right)} \to + \infty $ ;

2)$\left| {f_{S} \left( {t,x,s} \right)} \right| \le m_{1} + m_{2} \left|
{s} \right|^{p - 1},\,\,\,\,\,\forall t \in \left[ {a,b} \right],\,\,\,\,\,x
\in R^{n}; \quad 
m_1 ,m_2  \ge 0. $

\bigskip

3) there is a non-negative function\textit{} $a\left( {t,x,p}
\right)$\textit{,} $i \in I_{0} $ such, that if $\left\| {x}
\right\|_{W_{p}^{1} \left( {a,b} \right)} \le const$ then $\left\| {a_{i}
\left( { \cdot ,x\left( { \cdot}  \right),\dot {x}\left( { \cdot}  \right)}
\right)} \right\|_{L_{\infty}  \left( {a,b} \right)} \le const$ and the
function

\[
F(t,x,s) = f(t,x,s) + \sum\limits_{i \in I_0 } {a_i (t,x,s)f_{S_i
}^2 (t,x,s)}
\]

\noindent
is convex relatively\textbf{} $s$in\textbf{} $R^{n}$\textbf{.}

The following is proved.

\textbf{Theorem 1.} There exists a solution of the problem (\ref{eq1}) - (\ref{eq2}) under
the conditions 1) - 3).

The example of function\textit{} $f\left( {t,x,p} \right)$\textit{
}satisfying all conditions of the theorem 1 which is however not
convex relatively $s$ on $R^n$ is given.

The following multivariable problem also is considered

\begin{equation}
J(x) = \int\limits_D {f(t,\dot x(t)dt \to \min }
\end{equation}

Here

\[
D \subset R^{n}
t = \left( {t_{1} ,t_{2} ,...,t_{n}}  \right) \subset D
\dot {x}\left( {t} \right) = \left( {x_{t_{1}}  \left( {t} \right),x_{t_{2}
} \left( {t} \right),...,x_{t_{n}}  \left( {t} \right)} \right),
\quad
f:D \times R^n  \to R
\]

\noindent
is continuous on set of variables with the partial derivatives $f_{S_{1}}
,f_{S_{2}}  ,...,f_{S_{n}}  $.

The theorem of existence is formulated within similar conditions $\left(
{I_{0} = \{ 1,2,...,n\}}  \right)$.

Now we shall consider questions of resolvability of the optimal control
problem for the abstract evolutionary equations of the first order ([18]).

Let$H$ be Hilbert space, A$_{0} $ be linear closed
operator$
A_0 :D(A_0 ) \to H $, $\overline {D\left( {A_{0}}  \right)} = H$,
A (t) be the self-adjoint linear operator at every $t \in \left[
{0,T} \right]$, $A\left( {t} \right):D\left( {A_{0}} \right) \to
H$ and, satisfying inequalities

\[
m\left\| {A_{0} \varphi}  \right\|_{H} \le \left\| {A\left( {t}
\right)\varphi}  \right\|_{H} \le M\left\| {A_{0} \varphi}  \right\|_{H}
\varphi  \in D(A_0 ), m,n > 0.
\]

Through $H_{0} = D\left( {A_{0}}  \right)$ we shall designate Hilbert space
with norm

\[
\left\| {\varphi}  \right\|_{H_{0}} ^{2} = \left\| {\varphi}
\right\|_{H}^{2} + \left\| {A_{0} \varphi}  \right\|_{H}^{2} ,
\varphi  \in D(A_0 ).
\]

Let C be the linear continuous operator $C:H_{0} \to H_{1} $, $H_{1} $ is
Hilbert space , $V_{0} \subset U$ be the closed convex bounded set in
Hilbert space\textit{ U,}

\[
V = \left\{ {v:v \in L_{2} \left( {0,T,U} \right),v\left( {t} \right) \in
V_{0} ,\mathop {\forall} \limits^{0} t \in \left( {o,T} \right)} \right\},
\]

Operator\textit{ B (t, v, u)} for every $
t \in [0,T] $ and $v \in V_{0} $ continuously operates from $H_{0}
$ in H; operator\textit{ K (t, v, u, p)} for every $t \in \left[
{0,T} \right]$ continuously operates from $V_{0} \times H \times
H_{1} $ in R and at each

\[
(v,u,p) \in V \times L_2 (0,T;H) \times L_2 (0,T;H_1 ) \quad
K\left( {t,v\left( {t} \right),u\left( {t} \right),p\left( {t}
\right)} \right) \in L_{2} \left( {0,T} \right),
\]

Ô (u) is continuous functionals, defined in H.

The Banach space of functions u=u (t) is designated by W, belonging to
$L_{2} \left( {0,T;H_{0}}  \right)$, strict continuous on t in norm H,
having final norm

\[
\left\| {u} \right\|_{W} = \mathop {max}\limits_{0 \le t \le T} \left\|
{u\left( {t} \right)} \right\|_{H} + \left\| {A_{0} u} \right\|_{L_{2}
\left( {0,T;H} \right)} .
\]

Satisfy the subspace W is designated by$W_{0} $, for the elements which

\[
\int\limits_0^{T - h} {h^{ - 1} \left\| {u(t + h) - u(t)}
\right\|_H^2 dt \to 0} \textit{at}, h \to 0.
\]

The minimization problem for the functional is considered

\begin{equation}
\label{eq3}
J\left( {v} \right) = \int\limits_{0}^{T} {K\left( {t,v\left( {t}
\right),u\left( {t} \right),Cu\left( {t} \right)} \right)dt} + \Phi \left(
{u\left( {T} \right)} \right),
\end{equation}

\noindent
on set $V$ at conditions

\begin{equation}
\label{eq4}
u_{t} + A^{2}\left( {t} \right)u + B\left( {t,v,u} \right) = 0,
\end{equation}

\noindent
where $\varphi \in H$. The solution of the problem (\ref{eq3}) - (6) is understood
as a function $u = u\left( {t} \right) \in W_{0} $ satisfying equality

\begin{equation}
\label{eq5}
\begin{array}{l}
 \left( {u\left( {t} \right),\eta \left( {t} \right)} \right) +
\int\limits_{0}^{t} {\left[ { - \left( {u,\eta _{t}}  \right) + \left(
{A\left( {\tau}  \right)u,A\left( {\tau}  \right)\eta}  \right) + \left(
{B\left( {\tau ,v,u} \right),\eta}  \right)} \right]d\tau +}  \\
 = \left( {\varphi ,\eta \left( {0} \right)} \right),\mathop {\quad \dot
{\forall} }\limits^{} t \in \left( {0,T} \right), \\
 \end{array}
\end{equation}

\noindent
for $\forall \eta = \eta \left( { \cdot}  \right) \in L_{2} \left(
{0,T;H_{0}}  \right)
\eta _t  \in L_2 (0,T;H) $, $\left( { \cdot , \cdot}  \right)$ is
scalar product in H. It's suggested, that the solution of the
reduced problem (5) - (6) exists only and satisfies to the
estimation

\bigskip

. $\left\| {u} \right\|_{W} \le N,\forall v \in V$

\bigskip

Let

\[
K\left( {t,v,u,p} \right) = K_{0} \left( {t,u,p} \right) + K_{1} \left(
{t,v,u} \right),
\]

\[
B\left( {t,v,u} \right) = B_{0} \left( {t,u} \right) + B_{1} \left( {t,v,u}
\right).
\]

Operators\textit{ K (t, v, u, p), Ô (u), B (t, v, u)}, $u,p \in H$
have Frechet derivatives $K_{u} ,K_{p} ,\Phi _{z} ,B_{} $; the
operators K$_{1} $ (t, v, u), B$_{1} $ (t, v, u) have continuous
Frechet derivatives on $v \in V$. All these derivatives satisfy to
Lipschits condition in $
V \times L_2 (0,T;H) \times L_2 (0,T;H). $

Introduced the Hamilton - Potryagen functional for the problem (\ref{eq3}) - (6)

\begin{equation}
\label{eq6}
H\left( {t,v,u,\psi}  \right) = \left( {B_{1} \left( {t,v,u} \right),\psi}
\right) - K_{1} \left( {t,v,u} \right),
\end{equation}

\noindent
where $\psi = \psi \left( {t} \right) \in W_{0} $ is the solution of the
adjoint problem for the (\ref{eq3})-(6)

\bigskip

\begin{equation}
\label{eq7}
\begin{array}{l}
 - \psi _{t} \left( {t} \right) + A^{2}\left( {t} \right)\psi + B_{u} \left(
{t,v,u} \right)\psi = K_{u} \left( {t,v,u,Cu} \right) + \\
 + C^{\ast} K_{P} \left( {t,v,u,Cu} \right) \\
 \end{array}
\end{equation}

\begin{equation}
\label{eq8}
\psi \left( {t} \right) = - \Phi _{u} \left( {u\left( {T} \right)} \right).
\end{equation}

Let the following conditions fulfill:

\bigskip

\textit{1)} $K_{1} \left( {t,v,u} \right),B_{1} \left( {t,v,u}
\right) \quad B_{0} \left( {t,u} \right)$ are satisfy Lipschitz
condition on $v \in U,u \in H$ and

\[
K\left( {t,v,u,p} \right) \ge g\left( {t} \right),\quad \Phi \left( {u}
\right) \ge \mu > - \infty ,\quad g \in L_{1} \left( {0,T} \right);
\]

2) For all $
v_1 ,v_2  \in V_0 ,u \in H,\mathop {\forall t}\limits^0  \in
(0,T), $

3)

\[
\begin{array}{l}
 \left\| {B_{1} \left( {t,\frac{{v_{1} + v_{2}} }{{2}},u} \right) -
\frac{{1}}{{2}}B_{1} \left( {t,v_{1} ,u} \right) - \frac{{1}}{{2}}B_{1}
\left( {t,v_{2} ,u} \right)} \right\|_{H} \le \\
 \le \frac{{\chi _{1}} }{{4}}\left\| {v_{2} - v_{1}}  \right\|_{U}^{2} \\
 \end{array}
\]

\[
\begin{array}{l}
 K_{1} \left( {t,\frac{{v_{1} + v_{2}} }{{2}},u} \right) -
\frac{{1}}{{2}}K_{1} \left( {t,v_{1} ,u} \right) - \frac{{1}}{{2}}K_{1}
\left( {t,v_{2} ,u} \right) \le \\
 \le - \frac{{\chi} }{{4}}\left\| {v_{2} - v_{1}}  \right\|_{U}^{2} , \\
 \end{array}
\]

Where $\chi _{1} \ge 0,\chi _{2} \ge 0$;

3) $\left\| {\psi \left( {t} \right)} \right\|_{H} \le q,\forall v \in
V,\mathop {\forall t}\limits^{0} \in \left( {0,T} \right);$

\bigskip

4) It is possible to choose strongly converging subsequence from the
sequence $
u = u_n (t),\psi _{}  = \psi _n (t) $\textit{ }in $L_{2} \left(
{0,T;H} \right)$.

Last condition is fulfilled, for example, if $H_{0} \subset H$ is compact
and

\[
\left\| {u} \right\|_{W} + \left\| {u_{t}}  \right\|_{L_{2} \left( {0,T;H}
\right)} \le const,\forall v \in V.
\]

\textbf{Theorem 2.} Let $\chi > q\chi _{1} $. Then there exists a
solution of the problem (\ref{eq3}) - (6) and we can choose
strongly converging subsequence to the solution in $L_2(0,T,U)$
from any sequence.

It turns out the consequence from the theorem 2 which covers a wide class of
nonlinear problems.

\textbf{Corollary 1.} Let B$_{1} $ (t, v, u) be linearly on v and K$_{1} $
(t, v, u) is strict convex on v on $V_{0} $. Then the statement of the
theorem 2 is fulfilled.

The received results are applied to one optimal control problem for the
parabolic equations.

Let $D \in R^{n}$ be bounded domain with enough smooth boundary
$\Gamma \quad \Omega = D \times \left( {0,T} \right),S = \Gamma
\times \left[ {0,T} \right]x = \left( {x_{1} ,x_{2} ,...,x_{n}}
\right) \in D$ , V$_{0} $ be some closed, bounded, convex set
in\textbf{} $R^{n}$.

Let the functional be minimized

\begin{equation}
\label{eq9}
\begin{array}{l}
 J(v) = \int\limits_\Omega  {[K_0 (x,t,u(x,t),u_x (x,t))}  +  \\
  + K_1 (x,t,v(x,t),u(x,t))]dxdt + \int\limits_D {\Phi (u(x,T))dx}  \\
 \end{array}
\end{equation}

\noindent
on set

\[
V = \{ v = v\left( {x,t} \right):v = \left( {v_{1} ,v_{2} ,...,v_{m}}
\right) \in L_{2}^{\left( {m} \right)} \left( {\Omega}  \right),v\left(
{x,t} \right) \in V_{0} ,\mathop {\forall} \limits^{0} \left( {x,t} \right)
\in \Omega \} ,
\]

\noindent
at conditions

\begin{equation}
\label{eq10}
u_{t} - \sum\limits_{i,j = 1}^{n} {\left( {a_{ij} \left( {x,t}
\right)u_{x_{i}} }  \right)_{x_{j}}  + a\left( {x,t,v,u} \right) + f\left(
{x,t,u,u_{x}}  \right) = 0,\mathop {\forall \left( {x,t} \right) \in \Omega
}\limits^{0}}  ,
\end{equation}

\begin{equation}
\label{eq11}
u\left( {x,0} \right) = \varphi \left( {x} \right),x \in D,
\end{equation}

. $\left. {u} \right|_{S} = 0$ (14)

Here $
a_{ij} (x,t),i,j = \overline {1,n} ,a(x,t,v,u),f(x,t,u,p) $,
$\varphi \left( {x} \right)$ are the given measurable on $\left(
{x,t} \right) \in \Omega $ functions , continuous on $v \in V_{0}
u \in R,p \in R^{n}$; $a_(ij)=a_{ji}$ and $a_{ij} $ satisfy to a
condition of uniform ellipticity ;\textit{ a (x, t, v, u), f (x,
t, u, p)} and their partial derivatives relatively \textit{v, u}
and\textit{ p} satisfy to Lipschitz condition on $\left( {v,u,p}
\right) \in V_{0} \times R \times R^{n}$; $\varphi \in L_{2}
\left( {D} \right).$

Let the solution of the problem (\ref{eq10}) - (14) $u = u\left( {x,t} \right) \in
\dot {V}_{2}^{1,1/2} \left( {\Omega}  \right)$ exists.

From corollary1 it turns out

\textbf{Theorem 3.} Let $a(x,t,v,u)$\textit{}  be linearly
relatively v and $K_{1} \left( {x,t,v,u} \right)$ strict convex
relatively v on\textit{} $V_{0} $\textit{.} Then the solution of
the problem (\ref{eq9}) - (14) exists and we can choose
subsequence from any minimizing subsequence strongly converging in
$L_{2}^{\left( {m} \right)} \left( {\Omega}  \right)$ to the
solution a.

In spite of the fact that\textit{ a (x, t, v, u)} is linear
relatively\textit{ v}, there is a nonlinear term in the equation \textit{ f
(x, t, v, u, u}$_{x} $\textit{)} in which to take a limit, generally
speaking, is impossible. It raises the importance of the Theorem 3\textit{.}

If\textit{ a (x, t, v, u}) is not linear relatively v then the following
conditions are put

1)
\[
\begin{array}{l}
 \left| {a(x,t,\frac{{v_1  + v_2 }}{2},u) - \frac{1}{2}a(x,t,v_1 ,u) - \frac{1}{2}a(x,t,v_2 ,u)} \right| \le  \\
  \le \frac{{\chi _1 }}{4}\left\| {v_2  - v_1 } \right\|_{R^m }^2 ,\chi _1  \ge 0; \\
 \end{array}
\]

\bigskip

2) $\begin{array}{l}
 K_{1} \left( {x,t,\frac{{v_{1} + v_{2}} }{{2}},u} \right) -
\frac{{1}}{{2}}K\left( {x,t,v_{1} ,u} \right) - \frac{{1}}{{2}}K_{1} \left(
{x,t,v_{2} ,u} \right) \le \\
 \le - \frac{{\chi} }{{4}}\left\| {v_{2} - v_{1}}  \right\|_{R^{m}}^{2}
,\quad \chi \ge 0; \\
 \end{array}$

3) $\left| {\psi \left( {x,t} \right)} \right| \le q,\quad q \ge 0,\mathop
{\quad \mathop {\forall} \limits^{0} \left( {x,t} \right)}\limits^{} \in
\Omega ,\quad \forall v \in V.$

\bigskip

It turns out the resolvability of a problem of optimal control (\ref{eq9}) - (14)
at $\chi > q\chi _{1} $.

The optimal control problem for the abstract evolutionary equations of the
second order is investigated analogically.

The optimal control problem for Goursa-Darbou system also has been
considered ([14,18]).

Let the functional be minimized

\begin{equation}
\label{eq12}
\begin{array}{l}
 J\left( {v} \right) = \int\limits_{0}^{T} {\int\limits_{0}^{l} {K\left(
{x,t,u\left( {x,t} \right),u_{x} \left( {x,t} \right),u_{t} \left( {x,t}
\right),v\left( {x,t} \right)} \right)dxdt} +}  \\
 + \Phi \left( {u\left( {l,T} \right)} \right) \\
 \end{array}
\end{equation}

\noindent
at conditions

\begin{equation}
\label{eq13}
u_{xt} (x,t) = f(x,t,u(x,t),u_x (x,t),u_t (x,t),v(x,t)),(x,t) \in
Q ,
\end{equation}

\begin{equation}
\label{eq14}
u\left( {0,t} \right) = \varphi _{0} \left( {t} \right),t \in \left[ {0,T}
\right];u\left( {x,0} \right) = \varphi _{1} \left( {x} \right),\;\;x \in
\left[ {0,l} \right],
\end{equation}

Where\textit{ l, T} is the given positive numbers $\varphi _{0} \left( {0}
\right) = \varphi _{1} \left( {0} \right)$,

\[
Q = \{ \left( {x,t} \right):0 \le x \le l,\,\,\,0 \le t \le T\} f
= \left( {f^{1},f^{2},...,f^{n}} \right),\quad
\]
\[
\varphi _{j} = \left( {\varphi _{j}^{1} ,\varphi _{j}^{2}
,...,\varphi _{j}^{n}} \right),\quad j = 0,1, v = v(x,t) \in V ,
\]

\[
V = \{ v = \left( {v^{1},v^{2},...,v^{r}} \right):v \in L_{2}^{\left( {r}
\right)} \left( {Q} \right),v\left( {x,t} \right) \in D,\quad \dot {\forall
}\left( {x,t} \right) \in Q\} ,\quad D \subset E^{r}
\]

\noindent
is the convex closed bounded set.

Let's understand the vector function $u = u\left( {x,t} \right)$as the
solution of a problem (\ref{eq13}) - (\ref{eq14}) appropriated to control $v \in V$ which
has generalized derivatives $u\left( {x,t} \right),u_{t} \left( {x,t}
\right),u_{xt} \left( {x,t} \right) \in L_{2}^{\left( {n} \right)} \left(
{Q} \right)$ and, satisfying the equation (\ref{eq13}) almost everywhere in Q and to
conditions (\ref{eq5}) in sense of equality of the appropriate traces
$
u(0, \cdot ),u( \cdot ,0) $.

It is supposed, that functions $K\left( {x,t,u,p,q,v} \right),f^{i}\left(
{x,t,u,p,q,v} \right),i = \overline {1,n} $, $\Phi \left( {u} \right)$ and
their partial derivatives on\textit{ u, p, q, v} are continuous on set of
arguments and satisfy to Lipschitz condition on (\textit{u, p, q, v);}

Let

\[
K\left( {x,t,u,p,q,v} \right) = K_{0} \left( {x,t,u,v} \right) + K_{1}
\left( {x,t,u,p,q} \right)
\]

.$f^{i}\left( {x,t,u,p,q,v} \right) = f_{0}^{i} \left( {x,t,u,v} \right) +
f_{1}^{i} \left( {x,t,u,p,q} \right),i = \overline {1,n} $

The function

\[
H(x,t,u,v,\psi ) =  - K_0 (x,t,u,v) + (f_0 (x,t,u,v),\psi )
\]

\noindent
is introduced, where $\psi = \psi \left( {x,t} \right) = \left( {\psi
^{1}\left( {x,t} \right),\psi ^{2}\left( {x,t} \right),...,\psi ^{n}\left(
{x,t} \right)} \right)$ is a solution of the conjugate system.

The following conditions are put:

Let $\mathop {\forall} \limits^{0} \left( {x,t} \right) \in Q,\quad v_{i}
\in D,\quad i = 1,2,\quad \lambda \in \left( {0,1} \right),\quad u \in R$ ,

\begin{equation}
\label{eq15}
\begin{array}{l}
 \left| {f_{0} \left( {x,t,u,\frac{{v_{1} + v_{2}} }{{2}}} \right) -
\frac{{1}}{{2}}\lambda f_{0} \left( {x,t,u,v_{1}}  \right) -
\frac{{1}}{{2}}f_{0} \left( {x,t,u,v_{2}}  \right)} \right|_{E^{n}} \le \\
 \le \frac{{\chi _{1}} }{{4}}\left| {v_{2} - v_{1}}  \right|_{^{r}}^{2} , \\
 \end{array}
\end{equation}

\begin{equation}
\label{eq16}
\begin{array}{l}
 K_{0} \left( {x,t,u,\frac{{v_{1} + v_{2}} }{{2}}} \right) -
\frac{{1}}{{2}}K_{0} \left( {x,t,u,v_{1}}  \right) - \frac{{1}}{{2}}K_{0}
\left( {x,t,u,v_{2}}  \right) \le \\
 \le - \frac{{\chi} }{{4}}\left| {v_{2} - v_{1}}  \right|_{}^{2} , \\
 \end{array}
\end{equation}

\begin{equation}
\label{eq17}
\left| {\psi (x,t)} \right| \le R  , \quad  \quad \mathop
{\forall} \limits^{0} \left( {x,t} \right) \in Q
\end{equation}

\textbf{Theorem 4.} The solution of the problem (\ref{eq12}) - (\ref{eq14}) at $\chi > \chi
_{1} R$ exists and we can choose strongly converging subsequence to the
solution in $L_{2}^{\left( {r} \right)} \left( {Q} \right)$from minimizing
sequence.

From the theorem 4 is obtained the following.

\textbf{Corollary 2.} Let $f_{^{0}}^{i} \left( {x,t,u,v} \right)$ be linear
relatively v and $
K_0 (x,t,u,v),\mathop \forall \limits^0 (x,t) \in Q,u \in R^n $ is
strict convex on v in D. Then the statement of the theorem 4 is
valid.

At last we shall note the results concerning to investigation of existence
of the solution and a sufficient condition of an optimality for the
parabolic equations ([15]).

Let the functional be minimized

\begin{equation}
\label{eq18}
J\left( {v} \right) = \int\limits_{\Omega}  {K\left( {x,t,u\left( {x,t}
\right),,v\left( {x,t} \right)} \right)} dxdt + \int\limits_{D} {\Phi \left(
{x,u\left( {x,T} \right)} \right)dx} ,
\end{equation}

\noindent
at conditions

\begin{equation}
\label{eq19}
u_{t} - \sum\limits_{i,j = 1}^{n} {\left( {a_{ij} \left( {x,t}
\right)u_{x_{i}} }  \right)_{x_{j}} }  = f\left( {x,t,v,u,u_{x}}
\right),\quad \left( {x,t} \right) \in \Omega ,
\end{equation}

\begin{equation}
\label{eq20}
u\left( {x,0} \right) = \varphi \left( {x} \right),x \in D
\end{equation}

\begin{equation}
\label{eq21}
u\left( {\xi ,t} \right) = g\left( {\xi ,t} \right),\left( {\xi ,t} \right)
\in S
\end{equation}

 $
a_{i_j } (x,t),f(x,t,v,u,p),K(x,t,v,u),\Phi (x,z) $ are the
measurable on $\left( {x,t} \right) \in \Omega ,$ and continuous
on $v \in V_{0} ,u \in R,p \in R^{n},z \in R.$ Besides $a_{i_{j}}
= a_{j_{i}}  ,\,\,i,j = \overline {1,n} $ are satisfy to a
condition of uniform ellipticity; $\varphi \in L_{2} \left( {D}
\right),\,\,g \in L_{2} \left( {S} \right).$ Let the only solution
of the reduced problem (\ref{eq19}) - (\ref{eq21}) exist from
$
u(x,t) \in V_2^{1,0} (\Omega ) $, satisfying identity

\[
\begin{array}{l}
 \int\limits_D {u(x,t)\eta (x,t)dx + \int\limits_0^t {\int\limits_D {[ - u_{} \eta _t  + } } } \sum\limits_{i,j = 1}^n {a_{ij} (x,t)u_{x_i } \eta _{x_j } }  -  \\
  - f(x,t,v,u,u_x )\eta ]dxdt = \int\limits_D {\varphi (x)\eta (x,0)dx} ,\mathop {\mathop \forall \limits^0 }\limits^{} t \in [0,T] \\
 \end{array}
\]

For any $
\eta  = \eta (x,t) \in \mathop W\limits^0 _2{^{1,1}} (\Omega ) $.

Let the following conditions fulfill

1) $f\left( {x,t,v,u,p} \right)$ almost for every $\left( {x,t} \right) \in
\Omega $ is convex relatively $
(v,u,p) \in V_0  \times R \times R^n $;

2) Function $u \to f\left( {x,t,v,u,p} \right)$ decreases on R , $\mathop
{\forall} \limits^{0} \left( {x,t} \right) \in \Omega ,v \in V_{0} ,p \in
R^{n}$;

3) Function $K\left( {x,t,v,u} \right)\;$ and\textit{ Ô (x, z)} for $\mathop
{\forall} \limits^{0} \left( {x,t} \right) \in \Omega $ are convex
relatively $
(v,u) \in V_0  \times R $ and $z \in R$ ;

4) Functions $u \to K\left( {x,t,v,u} \right),\quad z \to \Phi \left( {x,z}
\right)\;$ increase on R , $\mathop {\forall} \limits^{0} \left( {x,t}
\right) \in \Omega $, $v \in V_{0} $.

\textbf{Theorem 5.} The solution of the problem
(\ref{eq19})-(\ref{eq21}) is convex relatively\textit{} $v \in
V$\textit{,} almost for all $\left( {x,t} \right) \in \Omega $, at
conditions 1), 2).

Using this result, it is proved

\textbf{Theorem 6.} Let the conditions 1)-4) be satisfy. Then the solution
of the problem (\ref{eq18}) - (\ref{eq21}) exists.

Introduced the Hamilton - Potryagen functional for the problem (\ref{eq18})-(\ref{eq21})

\[
H\left( {x,t,v,u,\psi}  \right) = f\left( {x,t,v,u,u_{x}}  \right)\Psi -
K\left( {x,t,v,u} \right).
\]

\textbf{Theorem 7}. Let $u^{\ast} \left( {x,t} \right)$ and $\Psi ^{\ast
}\left( {x,t} \right)$ be the solution of the basic and adjoint problem at
$
v = v^* (x,t) \in V $. Then it is sufficiency and necessary the
fulfillment of condition for optimality of the control $v^{\ast}
\left( {x,t} \right)$

\bigskip

\begin{equation}
\label{eq22}
\begin{array}{l}
 H\left( {x,t,v^{\ast} \left( {x,t} \right),u^{\ast} \left( {x,t}
\right),\psi ^{\ast} \left( {x,t} \right)} \right) = \\
 = \mathop {max}\limits_{v \in V_{0}}  H\left( {x,t,v,u^{\ast} \left( {x,t}
\right),\psi ^{\ast} \left( {x,t} \right)} \right),\quad \forall \left(
{x,t} \right) \in \Omega \,\,. \\
 \end{array}
\end{equation}

\bigskip

\subsection*{Reference}

\bigskip

1. Filippov A.F\textbf{.} Proceedings of MSU, math., mech.(Russia), 1959,N:
1, p. 25-52.

1. Lions J.L.\textbf{} Optimal control problems for distributed parameter
systems. Moscow, ``Nauka', 1972, 461 p.

3.Chesari L, Surjanarajany M.B.\textbf{} Journ. Optim.theory and
Appl.1970,v.31, N:3, p.307-415.

4. I.Ekland, R.Temam.\textbf{} Convex analysis and variational problems.
Amsterdam, New York, 1976.

5. Tolstonogov A.A. Proceedings of ASR (Russia), math. 2000,v.64,

N; 4,p.163-182.

6. Tolstonogov A.A.\textbf{} Math. Proceedings.(Russia) 2001,v.192,N:9,
p.125-142.

7. Raytum U.E.\textbf{} Differential equation,(Russia), 1983, v.19
,p.1044-1047.

8. Morduhovich B.S. Approximations\textbf{} methods in optimization and
optimal control. Moscow,''Nauka'', 1988,359 p.

9.\textbf{} Neustadt L.W.\textbf{} Math. Anal. And Appl. 1963, v.7, N;1,
p.101-107\textbf{.}

10. Polyak B.T. Proceedings of MSU, math., mech.(Russia),1968,N: 1, p.
30-40.

11.Baranger J. J.Math.Pures et Appl., 1973,v.52, N:4, p.377-405.

12.\textbf{} Iskenderov A.D.,Niftiyev A.A. Doclady of Acad. of Sciences of
ASA (Azerbaijan),1986, v.42, N:5, p.7-10.

13. Guliyev H.F. Avtomat. and Mechanics,(Russia), 1996, N:1, p.180-185.

14. Niftiyev A.A. Proceedings of BSU (Azerbaijan), 2000, N 2, p. 117-123.

15. Niftiyev A.A. Proceeding of Nijny- Novgorod University (Russia), 2001, N
23, p. 212-216.

16. Niftiyev A.A. Cybernetics and system analyses (Ukraina), 2001,N 6, p.
74- 79.

17. Niftiyev A.A., Gasimov Y.S. Control by boundaries and eigenvalue
problems with variable domains. Publ. House BSU , 2004, 185 p.

18. Niftiyev A.A. Existence of the solutions of non-convex optimization
problems, Publ. House BSU, 115 p.

19. Guliyev H.F., Niftiyev A.A. Proceeding of BSU 2004,¹ 3,p.21 - 28

\end{document}